\documentclass[12pt, twoside, leqno]{article}

\usepackage{amsmath,amsthm}
\usepackage{amssymb}

\usepackage{enumerate}

\newtheorem{thm}{Theorem}[section]
\newtheorem{cor}[thm]{Corollary}
\newtheorem{lem}[thm]{Lemma}
\newtheorem{alg}[thm]{Algorithm}

\theoremstyle{definition}
\newtheorem{defin}[thm]{Definition}
\newtheorem{rem}[thm]{Remark}

\newcommand{\Z}{\mathbb{Z}}
\newcommand{\N}{\mathbb{N}}

\DeclareMathOperator{\integ}{int}

\begin{document}

%%%%% To ease editing, for IMPAN journals add:

\baselineskip=17pt

%%%%%%%%%%%%%%%%

\title{Deterministic factorization of sums and differences of powers}

\author{Markus Hittmeir\\
University of Salzburg\\
Mathematics Department}

\date{}

\maketitle

%% Classification and key words; note that the 2010 classification is used:

\renewcommand{\thefootnote}{}

\footnote{The author is supported by the Austrian Science Fund (FWF): Project F5504-N26.}

\footnote{Address: Hellbrunnerstra{\ss}e 34, A-5020 Salzburg. E-Mail: markus.hittmeir@sbg.ac.at}

\footnote{2010 \emph{Mathematics Subject Classification}: 11A51.}

\footnote{\emph{Key words and phrases}: Factorization, Primes.}

\renewcommand{\thefootnote}{\arabic{footnote}}
\setcounter{footnote}{0}

%%%%%%

\vspace{-2.5cm}

\begin{abstract}
Let $a,b\in \N$ be fixed and coprime such that $a>b$, and let $N$ be any number of the form $a^n\pm b^n$, $n\in\N$. We will generalize a result of Bostan, Gaudry and Schost \cite{BosGauSch} and prove that we may compute the prime factorization of $N$ in 
\[
\mathcal{O}\Big{(}\textsf{M}_{\integ}\Big{(}N^{1/4}\sqrt{\log N}\Big{)}\Big{)},
\]
$\textsf{M}_{\integ}(k)$ denoting the cost for multiplying two $k$-bit integers. This result is better than the currently best known general bound for the runtime complexity for deterministic integer factorization.
\end{abstract}

\section{Introduction}
In \cite{CosHar} it has been proven that the deterministic and unconditional runtime complexity to compute the prime factorization of any natural number $N$ is in
\[
\mathcal{O}\Big{(}\textsf{M}_{\integ}\Big{(}\frac{N^{1/4}\log N}{\sqrt{\log\log N}}\Big{)}\Big{)},
\] 
where $\textsf{M}_{\integ}(k)$ denotes the cost for multiplying two $k$-bit integers. $\textsf{M}_{\integ}(k)$ can be bounded by $\mathcal{O}(k\log k 2^{\log^*k})$, where $\log^* k$ denotes the iterated logarithm (See \cite{Für}). The proof in \cite{CosHar} improves the well known approach of Strassen, which has been presented in \cite{Str}. The main idea of both methods is to use fast polynomial evaluation for computing parts of $\lfloor N^{1/2}\rfloor!$ to find a nontrivial factor of $N$.

In this paper, we will also apply this idea and combine it with a result of \cite{Hit} to improve the bound for numbers of certain shape, namely for sums and differences of powers. Our main theorem is the following:

{\thm{Let $a,b\in \N$ be fixed and coprime such that $a>b$, and define $P_{a,b}:=\{a^n\pm b^n: n\in\N\}$. Then, we may compute the prime factorization of any $N\in P_{a,b}$ in 
\[
\mathcal{O}\Big{(}\emph{\textsf{M}}_{\integ}\Big{(}N^{1/4}\sqrt{\log N}\Big{)}\Big{)}
\]
bit operations.}\label{c}}\\

We would like to point out that the theorem applies to some interesting subsets of $\N$, like Mersenne numbers or Fermat numbers.

\section{Preliminaries}
We briefly introduce the notions and results we will use in Section $3$ to prove Theorem \ref{c}. The following definitions describe the invertibility conditions required for fast polynomial evaluation provided by Theorem \ref{a}. They have been first introduced in \cite{BosGauSch}. Let $N\in\N$ and $\Z_N:=\Z/N\Z$.

\defin{Let $\alpha,\beta \in \Z_N$ and $d\in\N$. We say that $\textsf{h}(\alpha,\beta,d)$ is satisfied if the elements
\[
\beta, 2,...,d, (\alpha -d\beta),(\alpha-(d-1)\beta),...,(\alpha +d\beta)
\]
are invertible modulo $N$, and we define
\[
\textsf{d}(\alpha,\beta,d)=\beta2\cdots d(\alpha -d\beta)(\alpha-(d-1)\beta)\cdots(\alpha +d\beta).
\]
$\textsf{h}(\alpha,\beta,d)$ holds if and only if $\textsf{d}(\alpha,\beta,d)$ is invertible.
}

\defin{Let $\beta\in \Z_N$ and $e\in\N$. We say that $\textsf{H}(2^e,\beta)$ is satisfied if $\textsf{h}(2^{i},\beta,2^{i})$ and $\textsf{h}((2^{i}+1)\beta,\beta,2^{i})$ hold for each $0\leq i<e$. We define
\[
\textsf{D}(2^e,\beta)=\prod_{i=0}^{e-1}\textsf{d}(2^{i},\beta,2^{i})\textsf{d}((2^{i}+1)\beta,\beta,2^{i}).
\]
$\textsf{H}(2^{e},\beta)$ holds if and only if $\textsf{D}(2^e,\beta)$ is invertible.}

\begin{samepage}
{\lem{Let $f_0,...,f_{k-1}\in\Z_N$. Then we can decide if all $f_i$ are invertible modulo $N$ and, if not, find a noninvertible $f_i$ in 
\[
\mathcal{O}(k\emph{\textsf{M}}_{\integ}(\log N)+\log k \emph{\textsf{M}}_{\integ}(\log N)\log\log N)
\]
bit operations.}\label{d}}

\begin{proof}
See Lemma 12 in \cite{BosGauSch} for a proof.
\end{proof}
\end{samepage}
\pagebreak

Now let $H\in \Z_N[X]$ with $\deg H=1$. We define
\[
H_k(X)=H(X)H(X+1)\cdots H(X+k-1).
\]
Our main ingredients are the following two theorems.

{\thm{Let $\beta\in \Z_N$, $e\in\N$ and $k=2^e$. Assume that $\emph{\textsf{H}}(k,\beta)$ holds and that the inverse of $\emph{\textsf{D}}(k,\beta)$ is known. We may compute
\[
H_k(\beta),H_k(2\beta)...,H_k(k\beta)
\]
in $\mathcal{O}(\emph{\textsf{M}}_{\integ}(k\log (kN))+\emph{\textsf{M}}_{\integ}(\log N))$ bit operations.
}\label{a}}

\begin{proof}
Apply Proposition 7 in \cite{CosHar} with $\rho=1$.
\end{proof}

{\thm{Let $N\in\N$ be composite and $p$ a prime factor of $N$ with $p\leq b$ for some $b\leq N/5$. If $r,m\in\N$ such that $2\leq m< p$, $\gcd(N,m)=1$ and $r=p \mod m$, then the sets
\begin{align*}
&\{m^{-1}r-n \mod N:1\leq n\leq k\}\\
&\{-nk \mod N:1\leq n\leq k\}
\end{align*}
with $k=\lceil(b/m)^{1/2}\rceil$ are disjoint and there exist $i,j\in\{1,...,k\}$ such that $m^{-1}r-i\equiv -jk \mod p$.}\label{b}}

\begin{proof}
A proof can be found in \cite{Hit}, Corollary 4.4.
\end{proof}

{\cor{Let $N\in\N$ be composite and $p$ a prime factor of $N$ such that $p\leq b$ for some $b\leq N/5$. Let $r,m\in\N$ such that $2\leq m<p$, $\gcd(N,m)=1$ and $r=p \mod m$. Furthermore, define $H=X-m^{-1}r+1\in\Z_N[X]$ and set $k=\lceil (b/m)^{1/2} \rceil$. Then at least one of the elements 
\[
H_k(-k),H_k(-2k)...,H_k(-k^2)
\]
is noninvertible modulo $N$. Let $j\in\{1,...,k\}$ such that $H_k(-jk)$ is noninvertible, then $\gcd(-jk-m^{-1}r+i\mod N,N)$ yields a nontrivial factor of $N$ for some $i\in\{1,...,k\}$.
}\label{e}}

\begin{proof}
Theorem \ref{b} yields that there must exist $i,j\in\{1,...,k\}$ such that $m^{-1}r-i\equiv -jk \mod p$. This implies that the element $-jk-m^{-1}r+i$ is noninvertible modulo $N$. Hence, we conclude that the same holds for $H_k(-jk)=\prod_{l=1}^k -jk-m^{-1}r+l$. Since the sets in Theorem \ref{b} are disjoint, we get $p\mid\gcd(-jk-m^{-1}r+i\mod N,N)\neq N$, which proves the claim.
\end{proof}

{\lem{Let $N$ be a natural number  and $r,m\in\N$ with $m\geq 2$ such that $r=p\mod m$ for every prime divisor $p$ of $N$. Let $e\in\N$ such that $b:=4^em\leq N/5$. Knowing $r$ and $m$, one can compute a prime divisor $p$ of $N$ with $p\leq b$ or prove that no such divisor exists in
\[
\mathcal{O}(\emph{\textsf{M}}_{\integ}(2^{e}\log (N))+e\emph{\textsf{M}}_{\integ}(\log N)\log\log N)
\]
bit operations.}\label{f}}

\begin{proof}
If $m\geq q$ for the smallest prime factor $q$ of $N$, then $m=q$ or $r=q$. Assume $m<q$. Set $k:=2^{e}=\sqrt{b/m}$ and define $H=X-m^{-1}r+1\in\Z_N[X]$. We want to apply Theorem \ref{a} with $\beta=-k$ to compute the values
\[
H_k(-k),H_k(-2k)...,H_k(-k^2).
\]
In order to do this, we have to check if $\textsf{H}(k,-k)$ holds. It is easy to see that this the case if and only if $2,3,...,2^{e}+1$ and
\[
(2^{i}-2^{i}2^e),(2^{i}-(2^{i}-1)2^{e}),...,(2^{i}+2^{i}2^{e})
\]
are invertible modulo $N$ for each $0\leq i<e$. This list consists of $\mathcal{O}(k)$ easily computable elements in $\Z_N$, whose absolute values are bounded by $2^{e-1}+2^{e-1}2^{e}< 4^{e}<b\leq N/5$. Hence, they are all nonzero modulo $N$. By Lemma \ref{d}, we are able to decide if all of them are invertible modulo $N$ or, if not, find a noninvertible one in 
\[
\mathcal{O}(k\textsf{M}_{\integ}(\log N)+\log k \textsf{M}_{\integ}(\log N)\log\log N).
\]
Assume that we have found a noninvertible element, than we have also found a nontrivial factor of $N$ bounded by $k^2=4^{e}$. We are able to find a prime divisor of $N$ using trial division in $\mathcal{O}(k\textsf{M}_{\integ}(\log N))$ bit operations. In this case, the result is proven. Now assume that all of the elements above are noninvertible. We are able to compute $\textsf{D}(k,-k)\in\Z_N$ in $\mathcal{O}(k\textsf{M}_{\integ}(\log N))$ bit operations. The cost for computing its inverse are negligible. We now apply Theorem \ref{b} to compute $H_k(-nk)$ for $n=1,...,k$, and since $k<N$, we can do this in
\[
\mathcal{O}(\textsf{M}_{\integ}(k\log (N))+\textsf{M}_{\integ}(\log N))
\]
bit operations. Suppose that $N$ has a prime factor $p\leq b$ with $r=p\mod m$. Then by Corollary \ref{e}, there exists at least one $j\in\{1,...,k\}$ such that $H_k(-jk)$ is noninvertible modulo $N$. Using Lemma \ref{d}, we can find such an element in
\[
\mathcal{O}(k\textsf{M}_{\integ}(\log N)+\log k \textsf{M}_{\integ}(\log N)\log\log N).
\]
Let $H_k(-jk)$ be noninvertible modulo $N$, then Corollary \ref{e} yields that $\gcd(-jk-m^{-1}r+i\mod N,N)$ is nontrivial for some $i\in\{1,...,k\}$. Applying Lemma \ref{d} again, we are able to find such $-jk-m^{-1}r+i\mod N$ in
\[
\mathcal{O}(k\textsf{M}_{\integ}(\log N)+\log k \textsf{M}_{\integ}(\log N)\log\log N)
\]
bit operations. We know that $\gcd(-jk-m^{-1}r+i\mod N,N)$ is divisible by a prime divisor $p$ of $N$. This implies $-jk-m^{-1}r+i\equiv 0 \mod p$, hence $0\equiv mjk+r-mi \mod p$. We derive that $\gcd(mjk+r-mi,N)$ is nontrivial. The value $mjk+r-mi$ is bounded by $mk^2+r-m<mk^2=b$. Again, we use trial division to find a prime divisor of $N$. There are less than $k$ primes p smaller than $\lceil \sqrt{b}\rceil$ satisfying $r=p\mod m$, since they have to be of the form $mx+r$ for $x\in\{0,...,\lceil k/\sqrt{m}\rceil-1\}$. Therefore, the trial division can be done by $\mathcal{O}(k\textsf{M}_{\integ}(\log N))$ bit operations. This proves the claim.
\end{proof}

{\thm{Let $N\geq 400$ be a natural number  and $r,m\in\N$ with $m\geq 2$ such that $r=p\mod m$ for every prime divisor $p$ of $N$. Knowing $r$ and $m$, one can compute the prime factorization of $N$ in 
\[
\mathcal{O}\Big{(}\emph{\textsf{M}}_{\integ}\Big{(}\frac{N^{1/4}}{\sqrt{m}}\log N\Big{)}\Big{)}
\]
bit operations.}\label{g}}

\begin{proof}
We apply Lemma \ref{f} with $b=4^{e}m$ for $e\in\N$, starting with $e=1$. Lemma \ref{f} is applied with the same value of $b$ until no prime divisor of $N$ smaller than $b$ is found. Then we increase $e$ by $1$ and repeat. We do this until $b\geq \sqrt{N}$. Since $N\geq 400$, $b$ is always bounded by $4\sqrt{N}\leq N/5$.

If we run Lemma \ref{f} with value $b$, all prime divisors	 smaller than $b/4$ have already been detected. Since their product is bounded by $N$, we derive that the number of prime divisors between $b/4$ and $b$ and hence the number of runs of Lemma \ref{f} with the same value $b$ is bounded by $\mathcal{O}(\log N/\log b)$. The sum of all the terms of the form $e\textsf{M}_{\integ}(\log N)\log\log N)$ in the runtime complexity of Lemma \ref{f} is bounded by a polynomial in $\log N$ and hence negligible. We consider the sum of the other terms. Since $4^{e}m\geq \sqrt{N}$ implies $e\geq (\log N)/4-(\log m)/2$, we define $e_0:=\lceil(\log N)/4-(\log m)/2\rceil$ and get
\[
\sum_{i=1}^{e_0}\frac{\log N}{\log{(4^{i}m)}}\textsf{M}_{\integ}(2^{i}\log (N))
\leq \textsf{M}_{\integ}\Big{(}\log N \sum_{i=1}^{e_0}\Big{\lceil}\frac{\log N}{2i+\log m}\Big{\rceil} 2^{i}\Big{)}.
\]
Note that the inequality is a consequence of the facts that $k\leq \textsf{M}_{\integ}(k)$ and $\textsf{M}_{\integ}(k)+\textsf{M}_{\integ}(k')\leq \textsf{M}_{\integ}(k+k')$. 

Now we split the sum on the right side into $i\leq e_0/2$ and $i>e_0/2$. For $i\leq e_0/2$ we have $2^{i}\leq 2^{e_0/2}\in \mathcal{O}(N^{1/8}/m^{1/4})$, hence the first part of the sum is bounded by $\mathcal{O}((\log N)^2(N^{1/8}/m^{1/4}))$ and therefore negligible. We consider the main contribution by the summands with $i>e_0/2$. In these cases we have $2i+\log m>(\log N)/4-(\log m)/2+\log m>(\log N)/4$, hence the terms $\lceil \log N/(2i+\log m)\rceil$ are in $\mathcal{O}(1)$. We conclude that this part of the sum can be bounded by
\[
\mathcal{O}(2^{e_0})=\mathcal{O}(2^{\lceil(\log N)/4 -(\log m)/2\rceil})=\mathcal{O}(N^{1/4}/\sqrt{m}).
\]
which proves the claim.
\end{proof}

\rem{Let $N\in\N$ be odd. If we apply Lemma \ref{f} and Theorem \ref{g} with $m=2$ and $r=1$, we get the results of Lemma $13$ and Theorem $11$ in \cite{BosGauSch}.}

\section{Algorithm and Proof}
Let $a,b\in \N$ fixed and coprime such that $a>b$. We are interested in elements of $P_{a,b}^+:=\{a^n+b^n: n\in\N\}$ and $P_{a,b}^-:=\{a^n-b^n: n\in\N\}$. We define $P_{a,b}:=P_{a,b}^+\cup P_{a,b}^-$ and consider the following algorithm:

{\alg{Let $N\in P_{a,b}$. We can write either $N=a^m+b^m\in P_{a,b}^+$ or $N=a^m-b^m\in P_{a,b}^-$ for some $m\in\N$. Let $v:=1$ and take the following steps to compute the prime factorization of $N$:
\begin{enumerate}
\item{Apply trial division to remove all prime factors smaller than $400$ from $N$. Denote the resulting number by $N_0$. If $N\in P_{a,b}^+$, then set $N_1=N_0$. If $N\in P_{a,b}^-$, apply trial division to compute the prime factorization of $a-b$, remove all prime factors dividing $a-b$ from $N_0$ and denote the resulting number by $N_1$.}
\item{Apply trial division to compute all divisors of $m$. If $N\in P_{a,b}^+$, define $\mathcal{D}:=\{2d:d\mid m\}$. If $N\in P_{a,b}^-$, define $\mathcal{D}:=\{d:d\mid m\}$. For $l\geq 2$, let $d_1<d_2<...<d_l$ be the ordered list of all elements in $\mathcal{D}$.}
\item{Set $j=v$.}
\item{Compute $G_j=\gcd((ab^{-1})^{d_j}-1 \mod N_j,N_j)$. If $G_j=1$, then set $N_{j+1}=N_j$. If $1<G_j\leq N$, apply Theorem \ref{g} for $m=d_j$ and $r=1$ to compute the prime factorization of $G_j$, remove all prime factors dividing $G_j$ from $N_j$ and denote the resulting number by $N_{j+1}$. If $N_{j+1}=1$, stop. If not, set $v=j+1$ and go to Step $3$.}
\end{enumerate}
}}

\begin{proof}[Proof of Theorem \ref{c}] 
We have to prove that the algorithm is correct and runs in $\mathcal{O}(\textsf{M}_{\integ}(N^{1/4}\sqrt{\log N}))$.

First note that $N\in P_{a,b}^+$ implies $(ab^{-1})^{2m}\equiv 1\mod N$ and $N\in P_{a,b}^-$ implies $(ab^{-1})^{m}\equiv 1\mod N$. Hence, we have $N_{l+1}=1$ in any case and the algorithm always terminates. 

To prove correctness, it remains to show that the conditions in Theorem \ref{g} are always satisfied. Let $i\in\{1,...,l\}$ be arbitrary. If $N\in P_{a,b}^-$, prime factors dividing $a-b$ have already been removed in Step $1$, hence $1<G_i$ implies $d_i\geq 2$. Furthermore, all prime factors smaller than $400$ have been removed in Step $1$, hence $1<G_i$ implies $G_i\geq 400$. Let $p$ be any prime factor of $G_i$. We now have to prove that $1= p\mod d_i$. If $N\in P_{a,b}^-$, then $(ab^{-1})^m\equiv 1 \mod p$. Hence, the order $o$ of the element $ab^{-1}$ modulo $p$ is a divisor of $m$. We derive $o\in\mathcal{D}$. If $N\in P_{a,b}^+$, then $(ab^{-1})^m\equiv -1 \mod p$ and $(ab^{-1})^{2m}\equiv 1 \mod p$. Hence, the order $o$ of the element $ab^{-1}$ modulo $p$ is of the form $2d$ for some divisor $d$ of $m$. Again, we derive $o\in\mathcal{D}$. Now since $p$ divides $G_i$, we deduce $(ab^{-1})^{d_i}\equiv 1 \mod p$. Furthermore, $p$ divides $N_i$ and therefore has not been removed as prime factor in the previous steps. But this implies $(ab^{-1})^{d_j}\not\equiv 1\mod p$ for $1\leq j<i$, and we conclude that $o=d_i$. Since the order of any element is a divisor of the group order $p-1$, we derive $p\equiv 1 \mod d_i$ and the claim follows.

We now consider the runtime of the algorithm. The cost for Step $1$ is in $\mathcal{O}(1)$ and negligible. We are left with the task to discuss the runtime of Step $2$ and Step $4$.

Step $2$: Note that $N\geq a^m-b^m\geq(b+1)^m-b^m\geq b^{m-1}$. This implies $m\leq\log_b N+1=\log N/\log b+1\in\mathcal{O}(\log N)$ for $b\neq 1$. For $b=1$, it is also easy to show that $m$ is bounded by $\mathcal{O}(\log N)$. Hence, the cost to compute all divisors of $m$ and to bring them into the right order can be bounded by $\mathcal{O}((\log N)^{1+\epsilon})$ and is negligible.

Step $4$: The cardinality of $\mathcal{D}$ can be bounded by $\mathcal{O}(\log N)$. The cost for computing the greatest common divisors is negligible. Assume the computational worst case, in which we have to apply Theorem \ref{g} for every $j\in\{1,...,l\}$. We consider $1\leq j< l$. Then we have $a^{d_j}\equiv b^{d_j}\mod G_j$, hence $G_j\mid a^{d_j}-b^{d_j}$. If $N\in P_{a,b}^-$, we can write $d_j=m/k$ for some $k\geq 2$, and we deduce that
\[
G_j\leq a^{d_j}-b^{d_j}= a^{m/k}-b^{m/k}\leq (a^m-b^m)^{1/k}\leq (a^m-b^m)^{1/2}= N^{1/2}.
\]
As a consequence, the runtime of all the applications of Theorem \ref{g} for $1\leq j<l$ can be bounded by $\mathcal{O}(\textsf{M}_{\integ}(N^{1/8}(\log N))\log N)$ and is negligible. If $N\in P_{a,b}^+$, then we have $d_j=2m/k$ for some $k\geq 2$. Note that $G_j>1$ implies $d_j\neq m$, since $a^m\equiv -b^m\not\equiv b^m\mod p$ for every prime factor $p$ of $N$. We deduce that $k>2$ and therefore
\[
G_j\leq a^{d_j}-b^{d_j}=a^{\frac{2m}{k}}-b^{\frac{2m}{k}}\leq (a^m-b^m)^{\frac{2}{k}}\leq (a^m-b^m)^{\frac{2}{3}}< (a^m+b^m)^{\frac{2}{3}}=N^{2/3}.
\]
We hence conclude that in this case the runtime of all the applications of Theorem \ref{g} for $1\leq j<l$ can be  bounded by $\mathcal{O}(\textsf{M}_{\integ}(N^{1/6}(\log N))\log N)$ and is negligible. 

We now consider the runtime of Theorem \ref{g} for $j=l$. First note that $N\leq a^m+b^m< 2a^m$ implies $m> \log_a N-\log_a 2= \log N/\log a-\log_a 2$, which is in $\mathcal{O}(\log N)$. Hence, $m$ and $2m$ are both lower bounded by $\mathcal{O}(\log N)$. Assume the computational worst case, in which we have $G_{l}=N$. Then, the runtime is in
\[
\mathcal{O}\Big{(}\textsf{M}_{\integ}\Big{(}\frac{N^{1/4}}{\sqrt{\log N}}\log N\Big{)}\Big{)}=\mathcal{O}(\textsf{M}_{\integ}(N^{1/4}\sqrt{\log N})).
\]
This proves the result.
\end{proof}

{\cor{Let $N\in\N$ be a Mersenne number or a Fermat number. Then $N\in P_{2,1}$ and we may compute the prime factorization of $N$ in $\mathcal{O}\Big{(}\textsf{M}_{\integ}\Big{(}N^{1/4}\sqrt{\log N}\Big{)}\Big{)}$.}}


\begin{thebibliography}{HD}

%% Use the widest label as parameter.
%% Reference items can be numbered or have labels of your choice, as below.

%% In IMPAN journals, only the title is italicized; boldface is not used.
%% Our software will add links to many articles; for this, enclosing volume numbers in { } is helpful
%% Do not give the issue number unless the issues are paginated separately.

%%%%%%%%%%% To ease editing, use normal size:

\normalsize
\baselineskip=17pt

%%%%%%%%%%%%%
\bibitem[CH14]{CosHar} E. Costa, D. Harvey,
\emph{Faster deterministic integer factorization}, Math. Comp. 83, Pages 339-345, 2014.

\bibitem[H15]{Hit} M. Hittmeir, \emph{Digit Polynomials and their application to integer factorization}, http://arxiv.org/abs/1501.03078.

\bibitem[BGS07]{BosGauSch} A. Bostan, P. Gaudry, \'{E}. Schost, \emph{Linear recurrences with polynomial coefficients and application to integer factorization and Cartier-Manin operator}, SIAM J. Comput., 36(6):1777-1806, 2007.

\bibitem[S77]{Str} V. Strassen,
\emph{Einige Resultate \"uber Berechnungskomplexit\"at},
Jahresbericht der Deutschen Mathematiker-Vereinigung, Pages 1-8, 1976/77.

\bibitem[F09]{Für} M. F\"urer,
\emph{Faster integer multiplication}, SIAM J. Comput. 39(3):979-1005, 2009.

\end{thebibliography}
\end{document}